\documentclass[12pt,a4paper]{amsart}
\usepackage{amsfonts, amssymb}
\usepackage[all]{xy}
\numberwithin{equation}{section}

\def\B{\mathcal B}

\newtheorem{thm}{Theorem}
\newtheorem{Cor}{Corollary}

\newtheorem{lemma}{Lemma}
\newtheorem{conj}{Conjecture}

\newcommand{\Z}{\ensuremath{\mathbb Z}}

\begin{document}
\title[]{Lattice points on circles, squares in arithmetic progressions  and sumsets of squares}
\author{Javier Cilleruelo and Andrew Granville}
\address{Departamento de Matem\'aticas. Universidad Aut\'onoma de Madrid. 28049 Madrid}
\email{franciscojavier.cilleruelo@uam.es}
\address{D{\'e}partment  de Math{\'e}matiques et Statistique,
Universit{\'e} de Montr{\'e}al, CP 6128 succ Centre-Ville,
Montr{\'e}al, QC  H3C 3J7, Canada}
\email{andrew@dms.umontreal.ca}
\thanks {\it 2000 Mathematics Subject Classification:\rm 11N36
}

\begin{abstract}
We discuss the relationship between various additive problems concerning squares.
\end{abstract}

\maketitle
\section{Squares in arithmetic progression} Let $\sigma(k)$ denote the maximum of the number of
squares in $a+b,\dots ,a+kb$ as we vary over positive integers $a$
and $b$. Erd\H os conjectured that $\sigma(k)=o(k)$ which Szemer\'
edi \cite{Sz} elegantly proved as follows:\ If there are more than
$\delta k$ squares amongst the integers $a+b,\dots ,a+kb$ (where
$k$ is sufficiently large) then there exists four indices $1\leq
i_1<i_2<i_3<i_4\leq k$ in arithmetic progression such that each
$a+i_jb$ is a square, by Szemer\' edi's theorem. But then the
$a+i_jb$ are four squares in arithmetic progression, contradicting
a result of Fermat. This result can be extended to any given field
$L$ which is a finite extension of the rational numbers:\ From
Faltings' theorem we know that  there are only finitely many six
term arithmetic progressions of squares in $L$, so from Szemer\'
edi's theorem we again deduce that there are $o_L(k)$ squares of
elements of $L$ in any $k$ term arithmetic progression of numbers
in $L$. (Xavier Xarles \cite{X} recently proved that are never six
squares in arithmetic progression in $\Z [\sqrt d]$ for any $d$.)
 \smallskip

In his seminal paper {\it Trigonometric series with gaps}
\cite{Ru} Rudin stated the following conjecture:
\begin{conj}\label{arithprog}
$\sigma(k)=O(k^{1/2})$.
\end{conj}

It may be that the most squares appear in the
arithmetic progression $49+24i, \ 1\leq i\leq k$ once $k\geq 8$
yielding that $\sigma(k)=\sqrt{ \frac 83 k} +O(1)$.
Conjecture \ref{arithprog} evidently implies the following slightly weaker version:

\begin{conj}\label{arithprog2} For any $\varepsilon >0$ we have
$\sigma(k)=O(k^{1/2+\varepsilon})$.
\end{conj}

 Bombieri, Granville and Pintz
\cite{BGP} proved that $ \sigma(k)=O(k^{2/3+o(1)}),$ and recently
Bombieri and Zannier \cite{BZ} have proved that
$\sigma(k)=O(k^{3/5+o(1)})$.

\section{Rudin's approach}

Let $e(\theta):=e^{2i\pi\theta}$ throughout. The following
well-known conjecture was discussed by Rudin (see the end of section 4.6 in \cite{Ru}):

\begin{conj}\label{Rudin2} For any $2\le p<4$ there exists a constant $C_p$ such that, for any trigonometric polynomial $f(\theta)=\sum_{k}a_ke(k^2\theta)$ we have
\begin{equation}\label{fp}
\| f\|_p\le C_p\| f\|_2.
\end{equation}
\end{conj}

\noindent Here, as usual, we define $\| f\|_p^p:= \int_0^1
|f(t)|^p dt$ for a trigonometric polynomial $f$. Conjecture
\ref{Rudin2} says that the set of squares is a  $\Lambda (p)$-{\sl
set}  for any $2\le p<4$, where $E$ is a $\Lambda(p)$-set if there
exists a constant $C_p$ such that (\ref{fp}) holds for any  $f$ of
the form  $f(\theta)=\sum_{n_k\in E}a_ke(n_k\theta)$ (a so-called
{\sl  $E-$polynomial}). By H\" older's inequality we have, for
$r<s<t$, \begin{equation}\label{holder} \| f\|_s^{s(t-r)} \  \le \
\| f\|_r^{r(t-s)} \ \| f\|_t^{t(s-r)} ;
\end{equation}
taking $r=2$ we see that if $E$ is a $\Lambda(t)$-set then it  is a $\Lambda(s)$-set for all $s\leq t$.

Let $r(n)$ denote the number of representations of $n$ as the sum of two squares (of positive integers).  Taking $f(\theta)=\sum_{1\leq k\leq x} e(k^2\theta)$, we deduce that $\|f\|_2^2=x$, whereas  $\|f\|_4^4=\sum_n \# \{ 1\leq k,\ell\leq x: \ n=k^2+\ell^2\}^2 \geq \sum_{n\leq x^2} r(n)^2 \asymp x^2\log x$;
so we see that (\ref{fp}) does not hold in general for $p=4$.
\smallskip

Conjecture \ref{Rudin2} has not been proved for any $p>2$, though
Rudin \cite{Ru} has proved the following theorem.
\begin{thm}
If $E$ is a $\Lambda(p)$-set, then any arithmetic progression of
$N$ terms contains $\ll N^{2/p}$ elements of $E$. In particular,
if Conjecture \ref{Rudin2} holds for $p$ then
$\sigma(k)=O(k^{2/p})$.
\end{thm}
\noindent{\sl Proof}: \ We use  Fej\' er's kernel $\kappa_N(\theta):=\sum_{|j|\leq N} (1-\frac {|j|}N ) \ e(j\theta)$. Note that $\| \kappa_N\|_1=1$ and  $\| \kappa_N\|_2^2 =\sum_{|j|\leq N} (1-\frac {|j|}N )^2\ll N$ so, by (\ref{holder}) with $r=1<s=q<t=2$ we have
$\| \kappa_N\|_q^{q}   \ll 1^{2-q}  N^{q-1}$ so that
$\| \kappa_N\|_q \ll N^{1/p}$ where $\frac 1q + \frac 1p =1$.

Suppose that $n_1,n_2\dots ,n_\sigma$ are the elements of $E$ which lie in the arithmetic progression $a+ib, 1\leq i\leq N$.
If $n_\ell=a+ib$ for some $i, 1\leq i\leq N$ then $n_\ell=a+mb+jb$ where $m=[(N+1)/2]$ and $|j|\leq N/2$;
and so $1-\frac {|j|}N\geq \frac 12$.
Therefore, for $g(\theta):=\sum_{1\leq \ell\leq \sigma} e( n_\ell\theta)$, we have
$$
\int_0^1 g(-\theta) e((a+bm)\theta) \kappa_N(b\theta) d\theta   \geq \frac \sigma2 .
$$
On the other hand, we have $\| g\|_p \leq C_p \| g\|_2\ll \sqrt{\sigma}$ since $E$ is a $\Lambda(p)$-set and $g$ is an $E$-polynomial. Therefore, by H\" older's inequality,
\begin{equation} \label{calcul}
\left| \int_0^1 g(-\theta) e((a+bm)\theta) \kappa_N(b\theta) d\theta \right|  \leq \| g\|_p \ \| \kappa_N\|_q \ll \sqrt{\sigma} N^{1/p}
\end{equation}
and the result follows by combining the last two displayed equations.

\bigskip

It is known that Conjecture \ref{Rudin2} is true for polynomials
$f(\theta)=\sum_{k\le N}e(k^2\theta)$ and Antonio C\'ordoba
\cite{Co} proved that Conjecture \ref{Rudin2} also holds for
polynomials $f(\theta)=\sum_{k\le N}a_ke(k^2\theta)$ when the
coefficients $a_k$ are positive real numbers and non-increasing.

\section{Sumsets of squares}   For a given finite set of integers $E$ let $f_E(\theta)=\sum_{k\in E} e(k\theta)$. Mei-Chu Chang \cite{Chang}  conjectured that for any $\epsilon >0$ we have
$$
\| f_E\|_4\ll_\epsilon \| f_E\|_2^{1+\epsilon}
$$
for any finite set of squares $E$. As $\| f_E\|_4^4=\sum_n
r_{E+E}^2(n)$ where $r_{E+E}(n)$ is the number of  representations
of $n$ as a sum of two elements of $E$, her conjecture is
equivalent to:

\begin{conj}[Mei-Chu Chang]\label{Mei2} For any $\epsilon >0$ we have that
\begin{equation}
\| f_E\|_4^4 = \sum_n r_{E+E}^2(n)\ll
|E|^{2+\varepsilon} = \| f_E\|_2^{4+2\varepsilon}
\end{equation}
for any finite set $E$ of squares.
\end{conj}

We saw above that $\sum_n r_{E+E}^2(n)\gg |E|^{2} \log |E|$ in the special case  $E=\{1^2,\dots ,k^2\}$, so
conjecture \ref{Mei2} is sharp, in the sense one cannot entirely remove the $\epsilon$.

Trivially we have
$$
\| f_E\|_4^4= \sum_n r_{E+E}^2(n)\le \max
r_{E+E}(n)\sum_nr_{E+E}(n)\le |E|\cdot|E|^{2}=|E|^3
$$
for any set $E$; it is surprisingly  difficult to improve this estimate when $E$ is a set of squares. The best result such result is due to Mei-Chu Chang \cite{Chang} who proved
that
$$
\sum_n r_{E+E}^2(n) \ll |E|^{3}/\log^{1/12}|E|
$$
for any set $E$ of squares. Assuming a major conjecture of
arithmetic geometry we can improve Chang's result, in a proof
reminiscent of that in \cite{BGP}:

\begin{thm}\label{BombLang}
Assume the Bombieri-Lang conjecture. Then
$$
\sum_n r_{E+E}^2(n) \ll |E|^{\frac {11}4}
$$
\end{thm}

\noindent{\sl Proof}:\ One consequence of \cite{CH} is that there exists an integer $B$, such that if the Bombieri-Lang conjecture is true then for {\sl any} polynomial $f(x)\in \mathbb Z[x]$ of degree  five  or six which does not have repeated roots, there are no more than $B$ rational numbers $m$ for which $f(m)$ is a square.
For any given set of five elements $a_1^2,\dots ,a_5^2\in E$
consider all integers $n$ for which there exist $b_1^2,\dots ,b_5^2\in E$ with $n=a_1^2+b_1^2=\dots =a_5^2+b_5^2$. Evidently
$f(n)=(b_1\dots b_5)^2$ where $f(x)=\prod_{i=1}^5 (x-a_i^2)$, and so there cannot we more than $B$ such integers $n$. Therefore,
$$
\sum_n \binom{r_{E+E}(n)}5  = \sum_n \# \{ a_1^2,\dots ,a_5^2\in E:\
\exists b_1^2,\dots ,b_5^2\in E, \ \text{with} \ n=a_i^2+b_i^2, \ i=1,\dots ,5 \}
$$
$$
= \sum_{ a_1^2,\dots ,a_5^2\in E} \#\{ n:\
\exists b_1^2,\dots ,b_5^2\in E, \ \text{with} \ n=a_i^2+b_i^2, \ i=1,\dots ,5 \} \leq B \binom{|E|}5 .
$$
We have $\sum_n r_{E+E}(n)=|E|^2$; and so
$\sum_n r_{E+E}(n)^5 \ll \sum_n \binom{r_{E+E}(n)}5 +  \sum_n r_{E+E}(n) \ll |E|^5$. Therefore, by Holder's inequality, we have
$$
\sum_n r_{E+E}^2(n) \leq \left( \sum_n r_{E+E}(n) \right)^{3/4} \ \left( \sum_n r_{E+E}^5(n) \right)^{1/4} \ll |E|^{11/4}.
$$
\bigskip


\bigskip
\begin{conj}[Ruzsa]\label{Ruzsa}
If $E$ is a finite set on squares then, for every $\epsilon>0$ we
have $$|E+E|\gg |E|^{2-\epsilon}.$$
\end{conj}

\begin{thm}\label{Mei2toRuzsa}
Conjecture \ref{Mei2} implies Conjecture \ref{Ruzsa} (with the
same $\varepsilon$).
\end{thm}

\noindent{\sl Proof}:\ By the Cauchy-Schwarz inequality we have
$$
|E|^4 = (\sum_n r_{E+E}(n))^2\leq |E+E| \cdot \sum_n r_{E+E}^2(n)
$$
and the result follows.
\bigskip

\begin{thm}\label{Ruzsatoarithprog2}
Conjecture \ref{Ruzsa}  implies Conjecture \ref{arithprog2} (
$\varepsilon \rightarrow \frac{\varepsilon}{4-2\varepsilon} $)
\end{thm}

\noindent{\sl Proof}:\ If $E$ is a set of squares which is a subset of an
arithmetic progression $P$ of length $k$ then $E+E\subset
P+P$. From conjecture \ref{Ruzsa} we deduce that
$$
|E|^{2-\varepsilon}\ll|E+E|\le |P+P|=2k-1
$$
and the result follows.
\bigskip

In particular, theorems \ref{BombLang}, \ref{Mei2toRuzsa} and
\ref{Ruzsatoarithprog2} show  that the Bombieri-Lang conjecture
implies  $\sigma(k)\ll k^{4/5}$, which is easy to obtain by
directly applying the Bombieri-Lang conjecture to our arithmetic
progression. To do better than this suppose that there are
$\sigma_{r,s}$ squares amongst $a+ib, 1\leq i\leq k$ which are
$\equiv r \pmod s$; that is the squares amongst $a+rb+jsb, 0\leq
j\leq [k/s]$. This gives rise to $\binom{\sigma_{r,s}}6$ rational
points on the set of curves $y^2=x\prod_{j=1}^5 (x+n_j)$ for
$0\leq n_1<n_2<\dots <n_5\leq [k/s]$. Summing over all $r \pmod s$
and all $s>\sigma/10$ we get
$$
 \frac{k^5}{\sigma^5} \gg \binom{10k/\sigma }5 \gg \sum_{s>\sigma/10\ } \sum_{\ r\pmod s} \binom{\sigma_{r,s}}6 \gg  \sum_{s>\sigma/10} s \binom{[\sigma/s]}6\gg \sigma^2
$$
and we obtain $\sigma(k)\ll k^{5/7}$. Anyway  this upper bound was
improved unconditionally in \cite{BGP} and \cite{BZ}.

\bigskip

An {\sl affine cube} of dimension $d$ in $\mathbb Z$ is a set of
integers $\{ b_0 + \sum_{i\in I} b_i:\ I\subset \{ 1,\dots, d\}
\}$ for non-zero integers $b_0,\dots ,b_d$. In \cite{So}, Solymosi
states

\begin{conj}[Solymosi]\label{Solymosi}
There exists an integer $d>0$ such that there is no affine cube of
dimension $d$ of distinct squares.
\end{conj}

This conjecture follows from the Bombieri-Lang conjecture for if
there were an affine cube of dimension $d$ then for any $x^2\in \{
b_0 + \sum_{i\in I} b_i:\ I\subset \{ 3,\dots, d\} \}$ we have
that $x^2+b_1, x^2+b_2, x^2+b_1+b_2$ are also squares, in which
case there are $\geq 2^{d-2}$ integers $x$ for which
$f(x)=(x^2+b_1)(x^2+b_2)(x^2+b_1+b_2)$ is also square; and so
$2^{d-2}\leq B$, as in the proof of theorem \ref{BombLang}.

In \cite{So}, Solymosi gives a beautiful proof that for any set of
real numbers $A$, if $|A+A|\ll_d |A|^{1+\frac 1{2^{d-1}-1}}$ then
$A$ contains many affine cubes of dimension $d$. Therefore we
deduce a weak version of Ruzsa's conjecture from Solymosi's
conjecture:

\begin{thm}
Conjecture \ref{Solymosi} implies that there exists $\delta>0$ for
which $|A+A|\gg |A|^{1+\delta}.$
\end{thm}

\bigskip

The Erd\H os-Szemer\'edi conjecture states that for any set of integers $A$ we have
$$
|A+A|+|A\diamond A| \gg_\epsilon |A|^{2-\epsilon} .
$$
In fact they gave a stronger version, reminiscent of the Balog-Szemer\'edi-Gowers theorem:

\begin{conj}[Erd\H os-Szemer\'edi]\label{ErSz}
If $A$ is a finite set on integers and $G\subset A\times A$ with $|G|\gg |A|^{1+\epsilon/2}$ then
\begin{equation} \label{ErSzCon}
|\{ a+b: (a,b)\in G\}| + |\{ ab: (a,b)\in G\}|  \gg_\epsilon |G|^{1-\epsilon}.
\end{equation}
\end{conj}

Mei-Chu Chang \cite{Chang}  proved that a little more than
Conjecture \ref{ErSz} implies Conjecture \ref{Mei2}:

\begin{thm} \label{Chang1}
If (\ref{ErSzCon}) holds whenever $|G|\geq \frac 12 |A|$  then Conjecture \ref{Mei2} holds.
\end{thm}

\noindent{\sl Proof}:\ Let $B$ be a set of $k$ non-negative integers and $E=\{ b^2:b\in B\}$. Define
$G_M:=\{ (a_+,a_-):\ \exists b,b'\in B \ \text{with} \ a_+= b+b',\ a_-= b-b',\ \text{and} \ b^2-b'^2\in M\}$
where $M\subset E-E$; and so $A_M:=\{ a_+,a_-:\ (a_+,a_-)\in G_M\} \subset (B+B)\cup (B-B)$. Therefore
$|A_M|\leq 2|G_M|$.

Since $\{ a+a': (a,a')\in G\}$ and $\{ a-a': (a,a')\in G\} $ are subsets of $\{ 2b: b\in B\}$, they have $\leq k$ elements; and $\{ aa': (a,a')\in G\}\subset M$. Therefore (\ref{ErSzCon}) implies that $|G_M|^{1-\epsilon} \ll_\epsilon |M|+k$.  Since, trivially, $|G_M|\leq k^2$ we have
$\sum_{m\in M} r_{E-E}(m) = |G_M| \ll k^{2\epsilon}(|M|+k)$.

Now let  $M$  be the set of integers $m$ for which $r_{E-E}(m) \geq k^{3\epsilon}$,  so that  $\sum_{m\in M} r_{E-E}(m)\geq k^{3\epsilon}|M|$ and hence $\sum_{m\in M} r_{E-E}(m)   \ll k^{1+2\epsilon}$ by combining the last two equations.
Therefore, as $r_{E-E}(m) \leq k$,
$$
\| f_E \|_4^4 = \sum_{m} r_{E-E}(m)^2 \leq  \sum_{m\in E-E} k^{6\epsilon} + k\sum_{m\in M} r_{E-E}(m) \ll
k^{2+6\epsilon}.
$$
\bigskip

She also proves a further, and stronger result along similar lines:

\begin{conj}[Mei-Chu Chang]\label{ErSz2}
If $A$ is a finite set of integers  and $G\subset A\times A$ then
\begin{equation} \label{ErSzCon2}
|\{ a+b: (a,b)\in G\}| \cdot |\{ a-b: (a,b)\in G\}|  \cdot |\{ ab: (a,b)\in G\}|  \gg_\epsilon |G|^{2-\epsilon}.
\end{equation}
\end{conj}

\begin{thm}  Conjecture \ref{ErSz2} holds if and only if Conjecture \ref{Mei2} holds.
\end{thm}

\noindent{\sl Proof}:\ Assume Conjecture \ref{ErSz2} and define
$B, A$ and $G_M$ as in the proof of theorem \ref{Chang1}, so that $(\sum_{m\in M} r_{E-E}(m))^2=|G_M|^2\ll k^{2+2\epsilon} |M|$.
We partition $E-E$ into the sets $M_j:=\{ m:\ 2^{j-1}\leq r_{E-E}(m)<2^j\}$ for $j=1,2,\dots ,J:=[\log (2k)/\log 2]$; then
 $(2^{j-1}|M_j|)^2\leq (\sum_{m\in M_j} r_{E-E}(m))^2 \ll k^{2+2\epsilon} |M_j|$ so that $\sum_{m\in M_j} r_{E-E}(m)^2\leq
2^{2j}|M_j|\ll k^{2+2\epsilon}$. Hence
 $$
\| f_E \|_4^4 = \sum_{m} r_{E-E}(m)^2 < \sum_j \sum_{m\in M_j} r_{E-E}(m)^2 \ll J  k^{2+2\epsilon}\ll k^{2+3\epsilon},
$$
as desired.

Now assume Conjecture \ref{Mei2} and let $G_n:=\{ (a,b)\in G: \ ab=n\}$. Then
$|G|^2=\left( \sum_n |G_n|\right)^2 \leq |\{ ab: (a,b)\in G\}| \cdot \sum_n |G_n|^2$, while
$$
\sum_n |G_n|^2=
 \int_0^1 \left| \sum_{(a,b)\in G} e( 4abt) \right|^2 dt
= \int_0^1 \left| \sum_{(a,b)\in G} e((a+b)^2t) e(-(a-b)^2t) \right|^2 dt
$$
which, letting $E_\pm:=\{ r^2:\ r=a\pm b,\ (a,b)\in G\}$, is
$$
\leq \int_0^1 \left| \sum_{r^2\in E_+} e(r^2t) \sum_{s^2\in E_-} e(-s^2t) \right|^2 dt
\leq \| f_{E_+} \|^2 \| f_{E_-} \|^2
$$
by the Cauchy-Schwarz inequality.
Now $\| f_{E_\pm} \|^2\ll |\{ a\pm b: (a,b)\in G\}| \cdot |G|^{2\epsilon}$ by Conjecture \ref{Mei2}, and our result follows by combining the above information.
\bigskip

\section{Solutions of a quadratic congruence in short intervals}  We begin with a connection between additive combinatorics and the Chinese Remainder Theorem.
Suppose that $n=rs$ with $(r,s)=1$; and that for given sets of
residues $\Omega(r) \subset \mathbb Z/r\mathbb Z$ and $\Omega(s)
\subset \mathbb Z/s\mathbb Z$ we have  $\Omega(n) \subset \mathbb
Z/n\mathbb Z$ given by $m\in \Omega(n)$ if and only if there
exists $u\in  \Omega(r)$ and $v\in \Omega(s)$ such  that $m\equiv
u \pmod r$ and $m\equiv v\pmod s$.  When $(r,n/r)=1$ consider the
map which embeds $\mathbb Z/r\mathbb Z\to \mathbb Z/n\mathbb Z$ by
taking $u \pmod r$ and replaces it by $U \pmod n$ for which
$U\equiv u \pmod r$ and $U\equiv 0 \pmod {n/r}$; we write
$\Omega(r,n)$ the image of $\Omega(r)$ under this map. The key
remark, which follows immediately from the definitions, is that
$$
\Omega(n) = \Omega(r,n) + \Omega(s,n) .
$$
Thus if $n=p_1^{e_1}\dots p_k^{e_k}$ where the primes $p_i$ are distinct then
$$
\Omega(n) = \Omega(p_1^{e_1},n) + \Omega(p_2^{e_2},n) + \dots + \Omega(p_k^{e_k},n) .
$$

Particularly interesting is   where  $\Omega_f(n)$ is the set of
solutions $m\pmod n$  to $f(m)\equiv 0 \pmod n$, for given
$f(x)\in \mathbb Z[x]$. We are mostly interested in when there are
many elements of $\Omega_f(n)$ in a short interval where $f$ has
degree two.  {\sl A priori} this seems unlikely since  the
elements of the $\Omega(r,n)$ are so well spread out, that is they
have a distance $\geq n/r$ between any pair of elements since they
are all divisible by $n/r$.

The next theorem involves the distribution of the elements of $\Omega(n)$ in the simplest non-trivial case, in which each $\Omega(p_j^{e_j})$ has just two elements, namely $\{0,1\}$, so that  $\Omega(n)$ is the set of solutions of $x(x-1)\equiv 0\pmod n$.

\begin{thm} Let $\Omega(n)$ be the set of solutions of $x(x-1)\equiv 0\pmod n$. Then
\begin{enumerate}
    \item   $\Omega(n)$ has an element in the interval $(1,n/k+1)$.
    \item For any $\varepsilon >0$ there exists $n=p_1\dots p_k$
    such that  $\Omega(n)\cap (1,(\frac 1k - \varepsilon) n] =\emptyset$.
    \item For any $\varepsilon >0$ there exists $n=p_1\dots
    p_k$ such that if $x\in \Omega(n)$ then $|x|<\varepsilon n$.

\end{enumerate}
\end{thm}

\noindent{\sl Proof}. Let $\Omega(p_j^{e_j},n)=\{0,x_j\}$ where
$x_j\equiv 1\pmod{p_j^{e_j}}$ and $x_j\equiv 0\pmod{p_i^{e_i}}$
for any $i\ne j$. Then $\Omega(n)=\{0,x_1\}+\cdots +\{0,x_k\}$.
Let $s_0=n$ and $s_r$ be the least positive residue of $x_1+\cdots +x_r \pmod n$  for $r=1,\dots ,k$ so that $s_k=1$.
By the pigeonhole principle, there exists $0\leq l<m\leq k$ such that $s_l$ and $s_m$ lie in the same interval $(jn/k,(j+1)n/k]$,
and so $|s_l-s_m|<n/k$ with $s_m-s_l\equiv x_{l+1}+\cdots +x_m \pmod n \in \Omega(n)$.  If $s_m-s_l>1$ then we are done.
Now $s_k-s_1\equiv 1 \pmod n$ but is not $=1$, so  $s_m-s_l\ne 0,1$. Thus we must consider when $s_m-s_l<0$. In this case
$x_1+\dots +x_l+x_{m+1}+\cdots +x_k \pmod n \in \Omega(n)$
and is $\equiv s_k-(s_m-s_l)\equiv 1-(s_m-s_l)$, and the result follows.

\smallskip

To prove (2) take $k-1$ primes $p_1,\dots ,p_{k-1}>k$, and integers $a_j=[p_j/k]$  for  $j=1,\dots
,k-1$. Let $P=p_1\cdots p_{k-1}$ and determine $r\pmod P$ by the Chinese Remainder Theorem satisfying
$ra_j (P/p_j)\equiv 1\pmod {p_j}$ for $j=1,\dots ,k-1$. Now let
$p_k$ be a prime $\equiv r\pmod P$, and let $a_k$ the least positive
integer satisfying $a_kP\equiv 1\pmod{p_k}$. Let
$n=p_1\cdots p_k$ so that $x_i=a_in/p_i$ for $i=1,\dots ,k$.
Now $n/k\geq x_i>  n/k-n/p_i>0$ for $i=1,\dots,k-1$
and so since $x_1+\dots +x_k\equiv 1 \pmod n$ with $1\leq x_k<n$
we deduce that $x_1+\dots +x_k=n+1$ and therefore
$1+n/k\leq x_k< 1+n/k+n\sum_{i=1}^{k-1} 1/p_i$. Now elements of $\Omega(n)$ are of the form $\sum_{i\in
I}x_i$ and we have
$|\sum_{i\in I}x_i-n|I|/k|\le 1+2n\sum_{i=1}^{k-1} 1/p_i$,
and this is $<\epsilon n$ provided each $p_i>2k/\epsilon$. Finally, since the cases $I=\emptyset$ and
$I=\{1,\dots ,k\}$ correspond to the cases $x=0$ and $x=1$
respectively, we have that any other element is greater than
$(1/k-\varepsilon)n$.

\smallskip

To prove (3) we mimic the proof of (2) but now choosing non-zero integers $a_j$ satisfying
$|\frac{a_j}{p_j}|<\frac{\varepsilon}{2k}$ for $j=1,\dots ,k-1.$ This implies that $|a_k/p_k|<\varepsilon /2$ and then $|\sum_{i\in I}x_i|<\varepsilon n$.

\bigskip

In the other direction, we give a lower bound for the length of intervals containing $k$ elements of $\Omega(n)$.

\begin{thm}\label{thmcong+}  Let integer $d\geq 2$ be given, and suppose that for each prime power $q$ we are given a set of residues
$\Omega(q)\subset (\mathbb Z/q\mathbb Z)$ which contains no more than
$d$ elements. Let $\Omega(n)$ be determined for all integers $n$ using
the Chinese Remainder Theorem, as described at the start of this
section. Then, for any $k\geq d$, there are no more than $k$ integers $x\in \Omega(n)$ in any interval of length $n^{\alpha_d(k)}$, where
$\alpha_d(k)=\frac{1-\varepsilon_d(k)}d>0$ with $0<\varepsilon_d(k)=\frac{d-1}{k} +O( \frac{d^2}{k^2})$.
\end{thm}

\noindent{\sl Proof}. Let $x_1,\dots ,x_{k+1}$ elements of
$\Omega(n)$ such that $x_1<\cdots <x_{k+1}<x_1+n^{\alpha_d(k)}$.
Let $q$ a prime power dividing $n$. Each $x_i$ belongs to one of
the $d$ classes $\pmod q$ in $\Omega(q)$. Write $r_1,\dots , r_d$
to denote the number of these $x_i$ belonging to each class. Then
$\prod_{1\le i<j\le k+1}(x_j-x_i)$ is a multiple of
$q^{\sum_{i=1}^d \binom{r_i}{2}}$. The minimum of
$\sum_{i=1}^d\binom{r_i}{2}$ under the restriction $\sum_ir_i=k+1$
is $d\binom{r}{2}+rs$ where $r,s$ are determined by $k+1=rd+s,\
0\le s<d$. Finally
$$n^{\alpha_d(k) \binom{k+1}{2}}>\prod_{1\le i<j\le k+1}(x_i-x_j)>n^{d\binom{r}{2}+rs}$$ and we get a contradiction, by taking
$\alpha_d(k)=(d\binom{r}{2}+rs)/{\binom{k+1}{2}}$.

\

The next theorem is an easy consequence of the proof above.

\begin{thm} \label{shortint} If $x_1<\dots <x_k$ are solutions to the equation
$x_i^2\equiv a\pmod b$, then $x_k-x_1>b^{\frac 12-\frac
1{2\ell}}$, where $\ell$ is the largest odd integer $\leq k$.
\end{thm}

\noindent{\sl First proof}.  For any maximal prime power $q$
dividing $b$,  $(a,q)$ must be an square so we can write
$x_i=y_i\prod_q(a,q)^{1/2}$ with $y_i^2\equiv a'\pmod{q'}$ where
 $q'=q/(a,q)$ and  $(a',q')=1$. Let $\Omega(q')$ be the solutions of
$y^2\equiv a'\pmod{q'}$. Now, since $(a',q')=1$ we have that
$|\Omega(q')|\le 2$ and we can apply theorem \ref{thmcong+} to
obtain that $$x_k-x_1=(y_k-y_1)\prod_q (a,q)^{1/2}\ge \left
(\prod_q q/(a,q)\right )^{\alpha_2(k-1)}{\prod_q(a,q)}^{1/2}\ge
(\prod_qq)^{\alpha_2(k-1)}.$$ Now, notice that
$\alpha_2(k-1)=1/2-1/(2l)$ where $l$ is the largest odd number
$\le k$.
\medskip

\noindent{\sl Second proof}. Write $x_j^2=a+r_jb$ where
$r_1=1<r_2<\dots <r_k$ (if necessary, by replacing $a$ in the
hypothesis by $x_1^2-b$). Consider the $k$-by-$k$ Vandermonde
matrix $V$ with $(i,j)$th entry $x_j^{i-1}$. The row with $i=1+2I$
has $j$th entry $(a+r_jb)^I$; by subtracting suitable multiples of
the rows $1+2\ell, \ell<I$, we obtain a matrix $V_1$ with the same
determinant where the $(2I+1,j)$ entry is now $(r_jb)^I$.
Similarly the row with $i=2I+2$ has $j$th entry $x_j(a+r_jb)^I$;
by subtracting suitable multiples of the rows $2+2\ell, \ell<I$,
we obtain a matrix $V_2$ with the same determinant where the
$(2I+2,j)$ entry is now $x_j(r_jb)^I$. Finally we arrive at a
matrix $W$ by dividing out $b^I$ from rows $2I+1$ and $2I+2$ for
all $I$. Then the determinant of $V$, which is $\prod_{1\leq
i<j\leq k} (x_j-x_i)$, equals $b^{[ (k-1)^2/4]}$ times the
determinant of $W$, which is also an integer, and the result
follows.

The advantage of this new proof is that if we can get non-trivial lower bounds
on the determinant of $W$ then we can improve Theorem \ref{shortint}. We
note that $W$ has $(2I+1,j)$ entry $r_j^I$, and
$(2I+2,j)$ entry $x_jr_j^I$.

\medskip

\noindent{\sl Remark}: Taking $k=\ell$ to be the smallest odd
integer $ \geq \frac{\log b}{\log 4}$, then we can split our
interval into two pieces to deduce from Theorem \ref{shortint} a
weak version of Conjecture  \ref{cong}:\ There  are no more than
$\frac{\log 4b}{\log 2}$ solutions $x$ to the equation $x^2\equiv
a\pmod b$ in any interval of length $b^{1/2}$. From this it
follows that the number of solutions $x$ to the equation
$x^2\equiv a\pmod b$ in any interval of length $L$ is
$$\ll 1+ \frac{\log L}{\log \left( 1+\frac{b^{1/2}}L\right)} .$$
This result, with `$1/2$' replaced by `$1/d$', was proved for the
roots of any degree $d$ polynomial mod $b$ by Konyagin and Steger
in \cite{KS}.
%

\bigskip

A slightly improvement on the theorem above would have interesting
consequences.

\begin{conj}\label{cong} There exists a constant $N$ such that there are no
more than $N$ solutions $0<x_1<x_2<\cdots <x_N<x_1+b^{1/2}$ to the
equation $x_i^2\equiv a\pmod b$, for any given $a$ and $b$.
\end{conj}

\begin{thm}
Conjecture \ref{cong} implies Conjecture \ref{arithprog}.
\end{thm}

\noindent{\sl Proof}. Suppose that there are $\ell \gg k^{1/2}$
squares amongst $a+b, a+2b,\dots, a+kb$, which we will denote
$x_1^2<x_2^2<\dots <x_\ell^2$. By conjecture \ref{cong} we have
$x_\ell -x_1\geq  [(\ell-1)/N] b^{1/2}$, whereas $(k-1)b \geq
(x_\ell +x_1)(x_\ell -x_1) \geq (x_\ell -x_1)^2$. Therefore
$[(\ell-1)/N]^2 \leq (k-1)$ which implies that $\ell\leq
N(1+\sqrt{k-1})$.
\bigskip

Conjecture \ref{cong} would follow easily from theorem
\ref{thmcong+} if we could get the exponent $1/2$ for some $k$,
instead of $1/2-\varepsilon_2(k)$. Conjecture \ref{cong} can be
strengthened and generalized as follows:

\begin{conj}\label{cong+}  Let integer $d\geq 1$ be given, and suppose that for each prime power $q$ we are given a set of residues 
$\Omega(q)\subset (\mathbb Z/q\mathbb Z)$ which contains no more than
$d$ elements. $\Omega(b)$ is determined for all  integers $b$ using
the Chinese Remainder Theorem, as described at the start of this
section. Then, for any  $\epsilon>0$ there exists a constant
$N(d,\epsilon)$ such that for any integer $b$ there are no more
than $N(d,\epsilon)$ integers $n,\ 0\le n< b^{1-\epsilon}$ with
$n\in \Omega(b)$.
\end{conj}

In theorem \ref{thmcong+} we proved such a result with the exponent 
`$1-\epsilon$' replaced by `$1/d-\epsilon$'. We strongly believe Conjecture
\ref{cong+} with `$1-\epsilon$' replaced by `$1/d$', analogous to
Conjecture \ref{cong}. In a 1995 email to the second author, Bjorn Poonen
asked Conjecture \ref{cong+} with `$1-\epsilon$' replaced by
`$1/2$' for $d=4$; his interest lies in the fact that this would
imply the uniform boundedness conjecture for rational preperiodic
points of quadratic polynomials (see \cite{Poo}).

Conjecture \ref{cong+} does not cover the case 
$\Omega_f(b)=\{ m \pmod b: \ f(m)\equiv 0 \pmod n\}$ for all monic polynomials 
$f$ of degree $d$ since, for example, the polynomial
$(x-a)^d\equiv 0 \pmod {p^{k}}$ has got $p^{k-\lceil k/d\rceil }$ 
solutions $\pmod{p^k}$, rather than $d$.  One may avoid this difficulty
by restricting attention to squarefree moduli
(as in a conjecture  posed by Croot \cite{Cro}); or, to be less restrictive, 
note that if $f(x)$  has more than $d$ solutions $\pmod {p^{k}}$ then
$f$ must have a repeated root mod $p$, so that $p$ divides the discriminant 
of $f$:

\begin{conj}\label{cong+2}   Fix integer $d\geq 2$.
For any  $\epsilon>0$ there exists a constant
$N(d,\epsilon)$ such that for any monic $f(x)\in \mathbb Z[x]$  there are 
no more
than $N(d,\epsilon)$ integers $n,\ 0\le n< b^{1-\epsilon}$, with
$f(n)\equiv 0 \pmod b$ for any integer $b$ such that if $p^2$ divides $b$ then
$p$ does not divide the discriminant of $f$.
\end{conj}

\section{Lattice points on circles}

\begin{conj}\label{latticeprog} There exists $\delta>0$ and
integer $m>0$ such that
if $a_i^2+b_i^2=n$ with $a_i,b_i>0$  and $a_i^2\equiv a_1^2\pmod
q$ for $i=1,\dots ,m$ then $q=O(n^{1-\delta})$.
\end{conj}

\begin{thm}
Conjecture \ref{latticeprog} implies Conjecture \ref{arithprog}
\end{thm}

\noindent{\sl Proof}. Suppose that $x_1^2<\cdots <x_r^2$ are
distinct squares belonging to the arithmetic progression
$a+b,a+2b,\dots , a+kb$ with $(a,b)=1$, where $r>\sqrt{8lk}$, with
$l$ sufficiently large $>m$. We may assume that $(a,b)=1$ and that
$b$ is even. There are $r^2$ sums $x_i^2+x_j^2$ each of which
takes one of the values $2a+2b,2a+3b,\dots,2a+2kb$, and so one of
these values, say $n$, is taken $\ge r^2/(2k-1)>4l$ times. So we
can write $n=r_j^2+s_j^2$ for $j=1,2,\dots ,4l$ for distinct pairs
$(r_j,s_j)$, and let $v_j=r_j+is_j$. Note that $n\equiv 2\pmod 8$.
Let $\Pi =\prod_{1\le i<j\le 4l}(v_j-v_i)\ne 0.$ We will prove
that $|\Pi |\ge b^{4\binom{l}{2}}(n/2)^{\binom{2l}{2}}$, by
considering the powers of the prime divisors of $b$ and $n$ which
divide $\Pi $. Note that $(n/2,b)=1$.

Suppose $p^e\| b$ where $p$ is a prime, and select $w\pmod{p^e}$
so that $w^2\equiv a\pmod{p^e}$. Note that each $r_j,s_j\equiv w
\text{ or }-w\pmod{p^e}$: We partition the $v_j$ into four subsets
$J_1,J_2,J_3,J_4$ depending on the value of
$(r_j\pmod{p^e},s_j\pmod{p^e})$. Note then that $p^e$ divides
$v_j-v_i$ if $v_i,v_j$ belong to the same subset, and so $p^e$ to
the power $\sum_i \binom{|J_i|}{2}>4\binom{l}{2}$ divides $\Pi $.

Now let $p$ be an odd prime with $p^e\| n$. If $p\equiv 3\pmod 4$
then $p^{e/2}$ must divide each $r_j$ and $s_j$ so that
then $p^{(e/2){\binom{4l}{2}}}$ divides $\Pi$. If $p\equiv 1\pmod 4$
let us suppose $\pi$ is a prime in $\mathbb Z[i]$ dividing $p$;
then $\pi^{e_j}{\bar \pi}^{e-e_j}$ divides $v_j$ for some $0\le
e_j\le e$. If $e_i\le e_j$ we deduce that $\pi^{e_i}{\bar
\pi}^{e-e_j}$ divides $v_j-v_i$. We now partition the values of
$j$ into sets $J_0,\dots, J_e$ depending on the value of $e_j$.
The power of $\pi$ dividing $\Pi$ is thus
$\sum_{i=0}^ei\sum_{g=i+1}^e|J_i||J_g|+\sum_{i=0}^e(e-i){\binom{|J_i|}{2}}$, and the power of $\bar \pi$ dividing $\Pi$ is thus
$\sum_{i=0}^g
(e-g)\sum_{i=0}^{g-1}|J_i||J_g|+\sum_{i=0}^e(e-i){\binom{|J_i|}{2}}$. It is easy to show that $\sum_{0\le i<g\le
e}(i+e-g)m_im_g+e\sum_{i=0}^e{\binom{m_i}{2}}$, under the conditions
that $\sum_i m_i$ is fixed and each $m_i\ge 0$, is minimized when
$m_0=m_e,\ m_1=\cdots =m_{e-1}=0$. Therefore the power of $\pi$
plus the power of $\bar \pi$ dividing $\Pi$ is $\ge 2e{\binom{2l}{2}}$.

Finally $|r_j-r_i|,|s_j-s_i|\le (x_r^2-x_1^2)/(x_r+x_1)\le
(k-1)b/(2\sqrt{a+b})$, and so $|v_j-v_i|^2\le
(k-1)^2b^2/(2(a+b))$, giving that $|\Pi|\le
(k^2b^2/(2(a+b)))^{(1/2){\binom{4l}{2}}}$. Putting these all
together, as well as that $n>2(a+b)$ gives that
$2^{2l-1}(a+b)^{3l-1}\le k^{4l-1}b^{3l}$. However this implies that
$n\le 2k(a+b)\le 2^{1/2} k^{5/2}b^{1+1/(3l-1)} \ll k^{5/2}n^{(1+1/(3l-1))(1-\delta)} \ll k^{5/2}n^{1-\delta/2} $, for $l$ sufficiently large; and therefore $a+b<n\ll k^{O(1)}$.

Let $u_1,\dots u_d$ be the distinct integers in $[1,b/2]$ for which
each $u_j^2\equiv a \pmod b$, so that $d\asymp 2^{\omega(b)}$,
by the Chinese Remainder Theorem, where
$\omega(b)$ denotes the number of prime factors of $b$. The number of $x_i\equiv u_j \pmod {b/2}$ is
$\leq 1+ ((a+kb)^{1/2}-a^{1/2})/(b/2) \leq 1+2(k/b)^{1/2}$; and thus $r\ll 2^{\omega(b)} + k^{1/2} 2^{\omega(b)}/b^{1/2}$. This is  $\ll k^{1/2}$ provided $\omega(b)\ll \log k$, which happens when $ b\ll k^{O(\log\log k)}$ by the prime number theorem.
The result follows.

\medskip

Here is a flowchart of the relationships between the conjectures above:

\

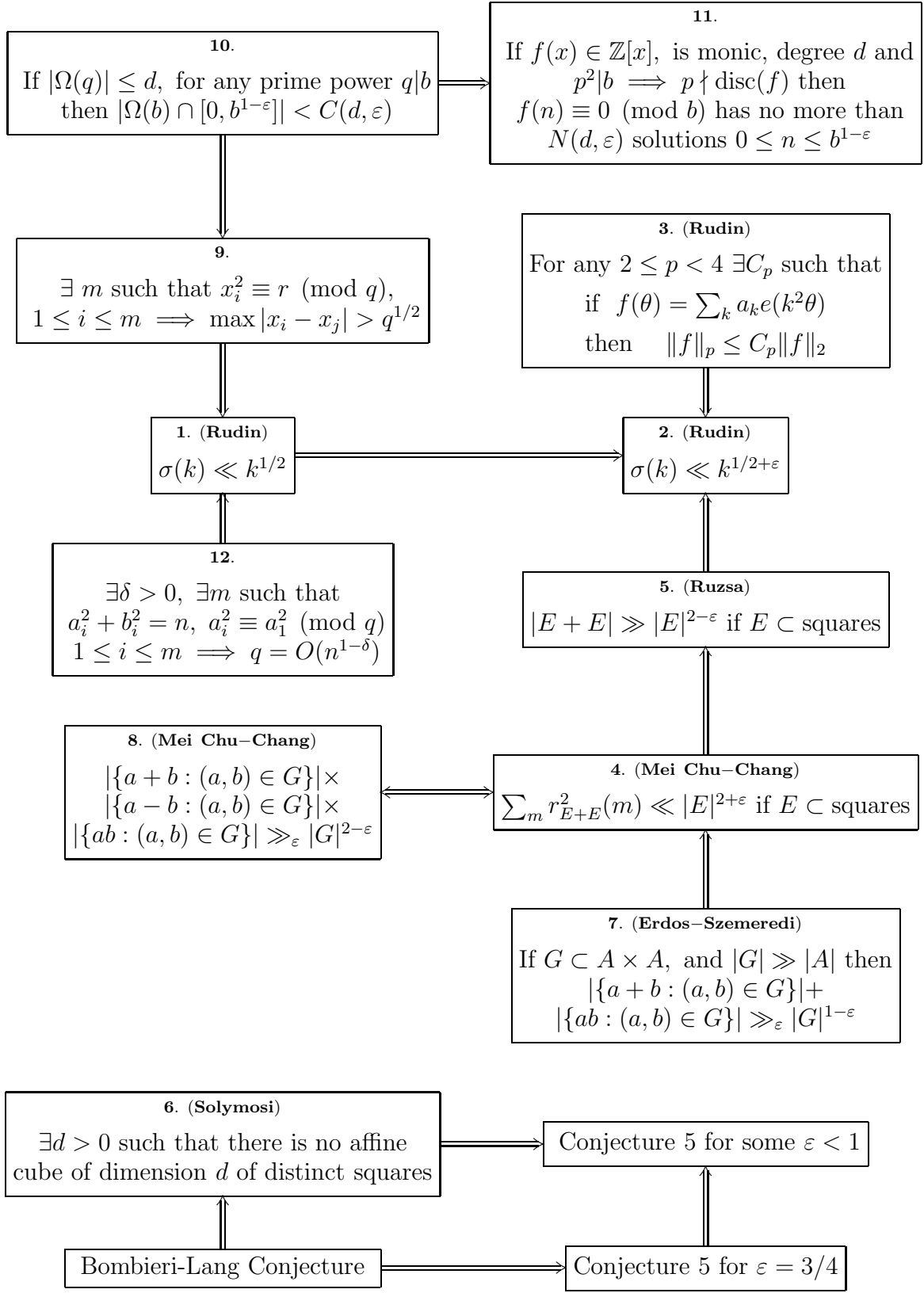
\bigskip

\noindent \xymatrix{ *+[F]{\ \begin{matrix} ^{\bf  10.\ }\\
\text{If } |\Omega(q)|\le d, \text{ for any prime power } q|b
\\\text{then } |\Omega(b)\cap [0,b^{1-\varepsilon}]|<C(d,\varepsilon)
\end{matrix}}\ar@2{=>}[d] & *+[F]{\ \begin{matrix} ^{\bf  11.\ }\\
\text{If } f(x)\in \Z[x], \text{ is monic, degree  } d \text{ and}
\\ p^2|b \implies p\nmid \text{disc}(f)\ \text{then} \\ f(n)\equiv 0\pmod b \text{ has no more than}\\ N(d,\varepsilon)
\text{ solutions } 0\le n\le b^{1-\varepsilon}
\end{matrix}}\ar@2{<=}[l]\\
*+[F]{\ \begin{matrix}
^{\bf  9.\ }\\
\exists \ m \text{ such that } x_i^2\equiv r\pmod q,\\ 1\le i\le m
\implies \max|x_i-x_j|>q^{1/2}\end{matrix}}\ar@2{=>}[d] & *+[F]{
\begin{matrix}
^{\bf 3.\ (Rudin)}\\ \text{For any }2 \le p<4\ \exists C_p \text{
such that}\medskip \\ \text{if }\ f(\theta)=\sum_ka_ke(k^2\theta)
\medskip\\\text{then }\quad \|f\|_p\le C_p \|f\|_2
\end{matrix}} \ar@2{=>}[d] &
\\*+[F]{
\begin{matrix}^{\bf  1.\ (Rudin)}\\ \sigma(k)\ll k^{1/2}\end{matrix}}\ar@2{=}[d]\ar@2{=>}[r]& *+[F]{
\begin{matrix}^{\bf  2.\ (Rudin)}\\  \sigma(k)\ll k^{1/2+\varepsilon}\end{matrix}} &\\
*+[F]{\ \begin{matrix} ^{\bf  12.\ }\\ \exists \delta>0,\ \exists
m
\text{ such that }\\ a_i^2+b_i^2=n,\ a_i^2\equiv a_1^2\pmod q \\
1\le i\le m \implies q=O(n^{1-\delta})
\end{matrix}}\ar@2{=>}[u] & *+[F]{\begin{matrix}
^{\bf  5.\ (Ruzsa)}\\ |E+E|\gg |E|^{2-\varepsilon} \text{ if }
E\subset \text{squares}\end{matrix}} \ar@2{=>}[u]\\*+[F]{
\begin{matrix}^{\bf  8.\ (Mei\ Chu-Chang)}\\ |\{a+b:(a,b)\in G\}|\times\\|\{a-b:(a,b)\in G\}|\times\\|\{ab:(a,b)\in G\}|\gg_{\varepsilon} |G|^{2-\varepsilon} \end{matrix} }
 \ar@2{<=>}[r]&
 *+[F]{\begin{matrix}
 ^{\bf 4.\ (Mei\ Chu-Chang)}\\ \sum_m r_{E+E}^2(m)\ll |E|^{2+\varepsilon} \text{ if } E\subset
 \text{squares}\end{matrix}}
 \ar@2{=>}[u] \\ & *+[F]{
 \begin{matrix}
 ^{\bf 7.\ (Erdos-Szemeredi)}\\ \text{If }G\subset A\times A,\text{ and } |G|\gg |A|\text{ then}\\ |\{a+b: (a,b)\in G\}|+\\ |\{ab:(a,b)\in G\}|\gg_{\varepsilon}|G|^{1-\varepsilon}
\end{matrix}
} \ar@2{=>}[u]\\ *+[F]{
\begin{matrix}
^{\bf  6.\ (Solymosi)}\\ \exists d>0 \text{ such that there is no
affine}\\ \text{cube of dimension } d \text{ of distinct
squares}\end{matrix} } \ar@2{=>}[r] & *+[F]{\begin{matrix}\text{
Conjecture } 5  \text{ for some } \varepsilon <1\end{matrix}
 }
 \\*+[F]{\text{ Bombieri-Lang Conjecture }
 }\ar@2{=>}[r] \ar@2{=>}[u] &
*+[F]{\begin{matrix}\text{Conjecture }5  \text{ for } \varepsilon
=3/4\end{matrix}
 } \ar@2{=>}[u]
 }

 \pagebreak

\begin{conj}\label{lattice1} For any $\alpha<1/2$, there exists a constant $C_{\alpha}$ such that for any $N$ we have
$$
\# \{ (a,b),\ a^2+b^2=n,\ N\le |b|<N+n^{\alpha}\}\le
C_{\alpha}.
$$
\end{conj}

A special case of interest is where $N=0$:\
\begin{equation}  \label{simplecase}
\# \{ (a,b),\ a^2+b^2=n,\  |b|<n^{\alpha}\}\le
C_{\alpha}.
\end{equation}
 Heath-Brown pointed out that one has to be careful in making an
analogous conjecture in higher dimension
as the following example shows:\  Select integer $r$ which has many
representations a the sum of two squares; for example, if $r$ is
the product of $k$ distinct primes that are $\equiv 1 \pmod 4$
then $r$ has $2^k$ such representations. Now let $N$ be an arbitrarily
large integer and consider the set of representations of $n=N^2+r$ as
the sum of three squares. Evidently we have $\geq 2^k$ such representations
in an interval whose size depends only on $k$, so is independent of
$n$.  However, one can get around this kind of example in formulating
the analogy to conjecture \ref{lattice1}  in 3-dimensions, since all
of these solutions live in a fixed hyperplane. Thus we may be able to
get a uniform bound on the number of such lattice points in
a small box, no more than
three of which live on the same hyperplane.

It is simple to prove (\ref{simplecase}) for any $\alpha\le 1/4$
(and Conjecture \ref{lattice1} for $\alpha\le 1/4$ with $N\ll
n^{\frac 12 -\alpha}$), but we cannot prove (\ref{simplecase}) for
any $\alpha >1/4$. Conjecture \ref{lattice1} and the special case
\ref{simplecase} are equivalent to the following conjectures
respectively:

\begin{conj}\label{lattice2} The number of lattice points $\{ (x,y)\in \mathbb Z^2:\ x^2+y^2=R^2\}$  in an arc of length
$R^{1-\epsilon}$ is bounded uniformly in $R$.
\end{conj}

\begin{conj}\label{lattice3} The number of lattice points $\{ (x,y)\in \mathbb Z^2:\ x^2+y^2=R^2\}$  in an arc of length
$R^{1-\epsilon}$ around the diagonal is bounded uniformly in $R$.
\end{conj}

Conjectures \ref{lattice1} and   \ref{lattice2} are simply a
rephrasing of one another, and obviously imply (\ref{simplecase})
and Conjecture \ref{lattice3}. In the other direction, if we have
points $\alpha_j:=x_j+iy_j$ on $x^2+y^2=R^2$ in an arc of length
$R^{1-\epsilon}$ then we have points
$\alpha_j\overline{\alpha_0}=a_j+ib_j$ satisfying
$a_j^2+b_j^2=R^2$  with $|b_j|\ll R^{1-\epsilon}$ contradicting
(\ref{simplecase}), and we have points
$(1+i)\alpha_j\overline{\alpha_0}$ on $x^2+y^2=2R^2$ in an arc of
length $\ll R^{1-\epsilon}$ around the diagonal, contradicting
Conjecture \ref{lattice3}.

The following result is proved in \cite{Ci-Co1}:

\begin{thm}\label{exponent}
There no more than $k$ lattice points $\{ (x,y)\in \mathbb Z^2:\
x^2+y^2=R^2\}$ in an arc of length $R^{\frac 12-\frac
1{4[k/2]+2}}$.
\end{thm}

\noindent{\sl Proof}. We may assume that $R^2=\prod_{p\equiv
1\pmod 4}p^e$, as the result for general $R^2$ is easily deduced from this case. Let $\mathbf p\overline{\mathbf p}$ be the
Gaussian factorization of $p$. Then each lattice point $\nu_i,\
1\le i\le k+1$ can be identified with a divisor of $R^2$ of the form $\nu_i=\prod_{\mathbf p} \mathbf p^{e_i}\overline{\mathbf
p}^{e-e_i}$. Therefore  $\nu_i-\nu_j$ is divisible by $\mathbf
p^{\min\{e_i,e_j\}}\overline{\mathbf p}^{\min\{e-e_i,e-e_j\}}$, so that
$|\nu_i-\nu_j|^2$ is divisible by $p^{e-|e_i-e_j|}$. Hence, since
$\sum_{1\le i<j\le k+1}|e_i-e_j|\le
e[\frac{k+1}2](k-[\frac{k+1}2])$, we have
$$
\prod_{1\le i<j\le k+1}|\nu_i-\nu_j|^2 \ge
\prod_pp^{\sum_{1\le i<j\le k+1}e-|e_i-e_j|}   \ge \left
(\prod_p p^e\right )^{{\binom{k+1}{2}}-[\frac{k+1}2](k-[\frac{k+1}2])}
$$
and the result follows.

\bigskip

It seems to be a difficult problem to decide whether the exponent $\frac
12-\frac 1{4[k/2]+2}$ is sharp for each $k$ in Theorem \ref{exponent}. We know that it is sharp for $k=1,2,3$ but we don't know what happens for larger $k$. More precisely:

\begin{enumerate}
    \item Obviously an  arc of length $\sqrt 2$ contains no more than one
lattice point; whereas the lattice points $(n,n+1),(n+1,n)$ lie on
an  arc of length  $\sqrt 2+o(1)$.
    \item It was shown in \cite{Ci3}  that an arc of length $(16R)^{1/3}$
contains no more than two lattice points. On the other hand the lattice
points $(4n^3-1,2n^2+2n), (4n^3,2n^2+1), (4n^3+1,2n^2-2n)$ lie on
an arc of length $(16R_n)^{1/3}+o(1)$.
    \item It was shown in  \cite{Ci-G} that an arc of length
    $(40+\frac{40}3\sqrt{10})^{1/3}R^{1/3}$, with $R>\sqrt{65}$,
contains no more than  three lattice points, whereas there exists an
infinite family of circles $x^2+y^2=R_n^2$ containing four lattice
points on an arc  of length \linebreak
$(40+\frac{40}3\sqrt{10})^{1/3}R_n^{1/3}+o(1)$. Other than in the examples arising from this
family, an arc of length
$(40+20\sqrt{5})^{1/3}R^{1/3}$ contains no more than three lattice points, whereas the four lattice points
$(x_0-2G_{n-2} ,y_0-2G_{n+1}), (x_0+G_{n-3},y_0+G_{n}), (x_0+G_{n-2},y_0+G_{n+1}), (x_0-G_{n-1}, y_0-G_{n+2})$,
where $x_0:=\frac 12 F_{3n+2},\ y_0=\frac 12 F_{3n-1},\ G_m=(-1)^mF_m$ and $F_m$ is the $m$th Fibonacci number, lie
on the circle $x^2+y^2=\frac 52F_{2n-2}F_{2n}F_{2n+2}=R_n^2$
on an arc of length  $(40+20\sqrt{5})^{1/3}R_n^{1/3}+o(1)$.
    \item Theorem \ref{exponent} is the best result known for all $k\ge 4$. In particular it implies that an arc of length
$R^{2/5}$ contains at most $4$ lattice points, and we do not know
whether the exponent $2/5$ can be improved: Are there infinitely many circles $x^2+y^2=R_n^2$ with four lattice points on an
arc  of length $\ll R_n^{2/5}$?
\end{enumerate}

\section{Incomplete trigonometric sums of squares}
The $L_4$ norm of a trigonometric polynomial has an interesting
number theory interpretation. For $f(\theta)=\sum_{n_k\in
E}a_ke(n_k\theta)$ we can write
\begin{displaymath}
\begin{split}
 \| f\|_4^4  &= \int_0^1 \left| \sum_{k}a_ke(n_k\theta) \right|^4 d\theta =\int_0^1\left|\sum_{m} \left(\sum_{n_k+n_j=m}a_ka_j
\right)e(m\theta)\right|^2 d\theta \\
   &=\sum_m \left|\sum_{n_k+n_j=m}a_ka_j\right|^2 \le \sum_m r_{E+E}(m)\sum_{n_k+n_j=m}|a_k|^2|a_j|^2 \\
   &\le \left(\sum_k |a_k|^2\right)^2 \max_mr_{E+E}(m)
\end{split}
\end{displaymath}
using the Cauchy-Schwarz inequality to obtain the first inequality, so that
\begin{equation} \label{ubound}
\| f\|_4 \leq \| f\|_2 \  \left(\sum_k \max_mr_{E+E}(m) \right)^{1/4} .
\end{equation}

If $E$ is the set of squares then $r_{E+E}(m)\leq \tau(m) \ll
m^{\varepsilon}$; so, by (\ref{ubound}), we have
$$
\| f\|_4\ll N^{\varepsilon}\| f\|_2
$$
for any $E$-polynomial $f$ where $E=\{1^2,\dots,N^2\}$.
Bourgain \cite{Bo} conjectured the more refined:

\begin{conj}\label{Bourgain} There exists a constant $\delta$ such that  for any $E$-polynomial $f$
where $E=\{1^2,\dots,N^2\}$, we have
$$\| f\|_4\ll \| f\|_2 (\log N)^{\delta} .$$
\end{conj}

\noindent Note that $\delta$ must be $\ge 1/4$; since we saw, in the second section, that
$\| f\|_4\sim C(\log N)^{1/4}\| f\|_2$ for $f(\theta)=\sum_{1\le k\le N}e(k^2\theta)$.

The corresponding conjecture when $f(\theta)=\sum_{k\in
E}e(k^2\theta)$ and $E\subset \{1^2,\dots ,N^2\}$ is the
following.

\begin{conj}\label{Bourgain'} There exists $C>0$ such that if
 $E\subset \{1^2,\dots ,N^2\}$ then $\sum_mr^2_{E+E}(m)\ll
|E|^2(\log N)^C$.
\end{conj}

Actually we can prove that both conjectures are equivalents.

\begin{thm}
Conjectures \ref{Bourgain} and \ref{Bourgain'} are equivalent.
\end{thm}

\noindent {\sl Proof:}  Conjecture \ref{Bourgain'} is a special case of Conjecture \ref{Bourgain}, so we must prove that Conjecture \ref{Bourgain} follows from Conjecture \ref{Bourgain'}. We may divide through the coefficients of $f$ by $\| f\|_2$ to ensure that $\| f\|_2=(\sum_{k}|a_k|^2)^{1/2}=1$, and therefore every $|a_k|\le 1$. Define $E_0=\{k,\ |a_k|\le N^{-1}\}$ and $E_j=\{ k,\ 2^{j-1}/N< |a_k|\le 2^{j}/N\}$ for all $j\geq 1$. Since $f=\sum_{j\geq 0} f_j$ (where each $f_j$ is the appropriate $E_j$-polynomial), we have $\|f \|_4\le
\sum_{j\geq 0} \|f_j\|_4$ by the triangle inequality.
By Conjecture \ref{Bourgain'} we have
$$
\|f_j\|_4^4=\sum_n \left|\sum_{\substack{k^2+j^2=n\\
k,j\in E_j}}a_ka_j\right|^2 \le (2^j/N)^4\sum_n r^2_{E_j+E_j}(n)\ll (\log N)^C(2^j/N)^4|E_j|^2.
$$

Now $\sum_{k\in E_j} |a_k|^2\asymp |E_j| (2^{2j}/N^2)$ for all $j\geq 1$, and $|E_0|/N^2,\ \sum_{k\in E_0} |a_k|^2 \ll 1/N$,
which imply that  $\sum_{j\geq 0} |E_j| (2^{2j}/N^2) \asymp 1$. Since $|E_j|=0$ for $j>\lceil \log_2 N\rceil$, we deduce that
$$
\frac 1{(\log N)^{C/4}} \sum_{j\geq 0} \|f_j\|_4  \ll
\sum_{j\geq 0} \frac{2^j|E_j|^{1/2}}N \ll
\left( \sum_{j=0}^{\lceil \log_2 N\rceil} 1   \sum_{j\geq 1} \frac{2^{2j}|E_j|}{N^2}   \right)^{1/2} \ll (\log N)^{1/2}.
$$
Therefore  Conjecture \ref{Bourgain} follows with $\delta=C/4+1/2$.

\medskip
Also we prove the following related result which slightly improves
on Theorem 2 of  \cite{Ci-Co1}.

\begin{thm} If $E=\{ k^2:\ N\le k\le N+\Delta\}$ with $\Delta\leq N$ and
$f(\theta)=\sum_{r\in E} e(r\theta)$, so that  $\| f\|_2^2\sim \Delta$,
then
$$\| f\|_4^4\asymp \Delta^2  + \Delta^3\cdot \frac{\log N}N .$$
In particular, $\| f\|_4\ll \| f\|_2$ if and only if $\Delta\ll (\log N)/N$.
\end{thm}

\noindent {\sl Proof}:\ Note that $\| f\|_2^2=|E|$ and
$$
\| f\|_4 = \sum_n r_{E+E}(n)^2 = 2|E|^2 -  |E| +
2\sum_n \left( \binom{r_{E+E}(n)}2 - \left[ \frac {r_{E+E}(n)}2
\right]\right) ;
$$
and that the sum counts twice the number of representations
$k_1^2+k_2^2=k_3^2+k_4^2$ with $N\le k_1,k_2,k_3,k_4\le N+\Delta$
and $\{ k_1,k_2\} \ne \{k_3,k_4\}$. Let $a+ib=$gcd$(k_1+ik_2,k_3+ik_4)$
and so
$k_1+ik_2=(a+ib)(x-iy)$ with $k_3+ik_4=(a+ib)(x+iy)u$ for some integers
$a,b,x,y$ where $u=1,-1,i$ or $-i$ is a unit.   Therefore
$k_1=ax+by, k_2=bx-ay$, and the four values of $u$
lead to the four possibilities $\{ k_3, k_4\} =
\{ \pm(bx+ay), \pm(ax-by)\}$. All four cases work much the same
so just consider $k_3=bx+ay, \ k_4=ax-by$.  Then
$N\leq ax=(k_1+k_4)/2,\ bx=(k_3+k_2)/2\le N+\Delta$
and $|by|=|k_1-k_4|/2,\ |ay|=|k_3-k_2|/2\le  \Delta/2$.
Multiplying through $a,b,x,y$ by $-1$ if necessary, we may assume $a>0$.
Therefore  $1+\Delta/N\geq b/a \geq (1+\Delta/N)^{-1}$ so that
$$
b=a+O(a\Delta/N),\ N/a\leq x\le N/a+\Delta/a,\ |y|\leq \Delta/2a .
$$
We may assume that  $a<\Delta$ else $y=0$ in which case
$\{ k_1,k_2\} \ne \{k_3,k_4\}$. Therefore, for a given $a$
the number of possibilities for $b,x$ and $y$ is
$\ll (a\Delta/N)(\Delta/a)^2=\Delta^3/aN$. Summing up over all $a,
1\leq a\leq \Delta$, gives  that $\| f\|_4\ll \Delta^3(\log \Delta)/N$.

On the other hand if integers $a,b,x,y$ satisfy
$$
a\in [7N/\Delta, \Delta/2],\ b\in [a(1-\Delta/7N), a], \
ax\in [N+\Delta/2,N+2\Delta/3],\ ay\in [1,\Delta/3],
$$
then $N\leq k_1=ax+by<k_2=bx-ay,\ k_3=bx+ay<k_4=ax-by\leq N+\Delta$
for $\Delta\leq N/3$, and so
$\| f\|_4\gg \Delta^2+ \Delta^3(\log (\Delta^2/N))/N$.
\medskip

\begin{conj}\label{Cico1}The exists $\eta$ such that for  any $E$-polynomial $f$ with $E=\{N^2,\dots ,(N+N/(\log N)^{\eta})^2\}$, we have
$$\| f\|_4\ll\| f\|_2.$$
\end{conj}

\medskip

Conjecture \ref{Cico1} probably holds with  $\eta=1$.
If  $E=\cup_{i=1}^r E_i$ then we can write any $E$-polynomial $f$ as $f=\sum_{i=1}^r f_i$, and  by  the triangle inequality we have $|f|^4\leq \sum_{i=1}^r |f_i|^4$. Therefore Conjecture \ref{Cico1} implies Bourgain's
Conjecture \ref{Bourgain} with $\delta=\eta/2$.

\medskip

In \cite{Ci-Co1} the following weaker conjecture was posed.

\begin{conj}\label{Cico2}For any $\alpha<1$, for any trigonometric polynomial $f$
with frequencies in the set $\{N^2,\dots ,(N+N^{\alpha})^2\}$, we
have
$$\| f\|_4\ll_{\alpha} \| f\|_2.$$
\end{conj}

Conjecture \ref{Cico2} is trivial for $\alpha\le 1/2$, yet is
completely open for any $\alpha >1/2$. From (\ref{ubound}) we
immediately deduce:

\begin{thm}
Conjecture \ref{lattice1}  implies Conjecture \ref{Cico2} .
\end{thm}

The next conjectures \ref{Cico1'} and \ref{Cico2'} correspond to
conjectures \ref{Cico1} and \ref{Cico2}, respectively, in the particular case
$f(\theta)=\sum_{k^2\in E}e(k^2\theta)$ and are also open.

\begin{conj}\label{Cico1'}
There exists $\delta >0$ such that if $E\subset \{k^2,\ N\le k\le
N+N/\log^{\delta}N\}$ then $\sum_m r^2_{E+E}(m)\ll |E|^2$.
\end{conj}

\begin{conj}\label{Cico2'}
If $E\subset \{k^2,\ N\le k\le N+N^{1-\varepsilon}\}$ then
$\sum_mr^2_{E+E}(m)\ll |E|^2$.
\end{conj}

We now give a flowchart describing the relationships between the
conjectures in the second half of the paper.
\bigskip

\noindent 
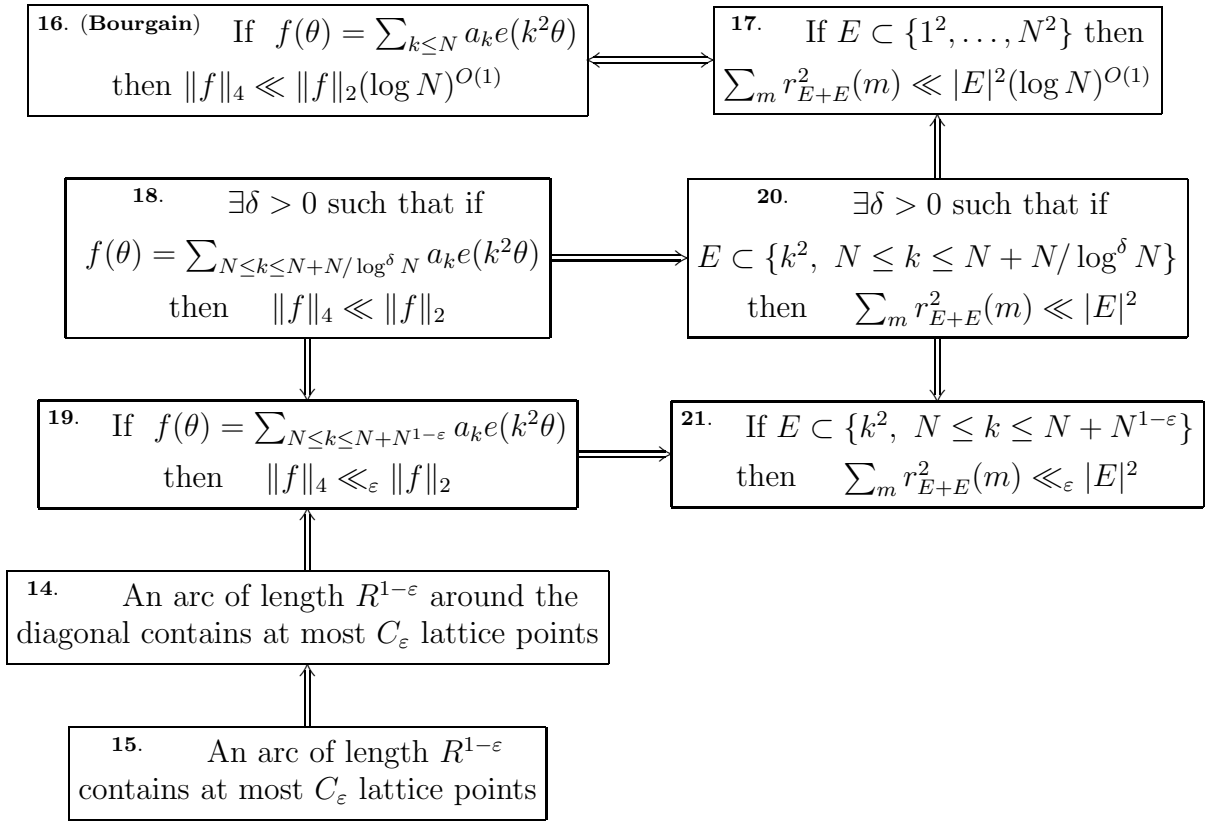
\xymatrix{
 *+[F]{\begin{matrix}^{\bf 16.\ (Bourgain)} \quad \text{If }\ f(\theta)=\sum_{k\le N}a_ke(k^2\theta)\medskip \\ \text{then }\| f\|_4\ll \|f\|_2(\log N)^{O(1)}\end{matrix}}
 \ar@2{<=>}[r]&
 *+[F]{\begin{matrix}^{\bf 17.} \quad \text{ If } E\subset \{1^2,\dots ,N^2\} \text{ then}\medskip \\ \sum_m r_{E+E}^2(m)\ll |E|^{2}(\log
 N)^{O(1)}\end{matrix}
 }\\
*+[F]{\begin{matrix}^{\bf 18.} \qquad \exists \delta >0 \text{ such that if}\medskip\\
\ f(\theta)=\sum_{N\le k\le N+N/\log^{\delta} N}a_ke(k^2\theta)\medskip \\
\text{then }\quad \| f\|_4\ll \|f\|_2  \end{matrix}
}\ar@2{=>}[r]&*+[F]{\begin{matrix}^{\bf 20.} \qquad \exists \delta
>0\text{  such that if }\medskip\\ E\subset \{k^2,\ N\le k\le
N+N/\log^{\delta}N\}
\medskip \\\text{ then }\quad \sum_m r_{E+E}^2(m)\ll |E|^{2}\end{matrix}
 }\ar@2{=>}[u]
\\
*+[F]{\begin{matrix}^{\bf 19.} \quad \text{If }
\ f(\theta)=\sum_{N\le k\le N+N^{1-\varepsilon}}a_ke(k^2\theta)\medskip \\
\text{then }\quad \| f\|_4\ll_{\varepsilon} \|f\|_2  \end{matrix}
}\ar@2{<=}[u]\ar@2{=>}[r] & *+[F]{\begin{matrix}^{\bf 21 .} \quad
\text{If } E\subset \{k^2,\ N\le k\le N+N^{1-\varepsilon}\}
\medskip \\
\text{then }\quad \sum_m r_{E+E}^2(m)\ll_{\varepsilon} |E|^{2} \end{matrix} }\ar@2{<=}[u]\\
*+[F]{\begin{matrix}^{\bf 14.}\qquad   \text{An arc of length } R^{1-\varepsilon} \text{ around the }\\
\text{diagonal contains at most  } C_{\varepsilon} \text{ lattice
points}
\end{matrix}}\ar@2{=>}[u]& \\
*+[F]{\begin{matrix}^{\bf 15.}\qquad \text{An arc of length } R^{1-\varepsilon}\\
\text{contains at most  } C_{\varepsilon} \text{ lattice points}
\end{matrix}}\ar@2{=>}[u]
 }

\bigskip

\section{Sidon sets of squares}

A set of integers $A$ is called a Sidon set if we have
$ \{a,b\}=\{c,d\}$ whenever $a+b=c+d$ with $a,b,c,d\in A$.
More generally $A$ is a $B_2[g]$-set if there are $\leq g$ solutions to $n=a+b$ with $a,b\in A$, for all integers $n$
(so that a Sidon set is a $B_2[1]$-set). The set of squares is
not a Sidon set, nor a $B_2[g]$-set for any $g$; however it is
close enough that this inspired Rudin in his
seminal article \cite{Ru}, as well as this paper.

One question is to find the largest Sidon set
$A\subset \{1^2,\dots,N^2\}$. Evidently $A=\{ (N-[\sqrt N]+k)^2,\ k=0,\dots ,[\sqrt N]-1\}$ is a Sidon
set of size $[\sqrt N]$. Alon and Erd\H os \cite{Alon} used the probabilistic method to obtain a Sidon set $A\subset \{1^2,\dots,
N^2\}$ with $|A|\gg_{\varepsilon} N^{2/3-\varepsilon}$ (and Lefmann and Thiele \cite{Le} improved this to $|A|\gg_{\varepsilon} N^{2/3}$).

We ``measure'' the size of infinite Sidon sets  $\{a_k\}$ by
giving an upper bound for $a_k$. Erd\H os and Renyi \cite{ER}
proved that there exists an infinite $B_2[g]$-set $\{ a_k\}$ with $a_k\ll k^{2+\frac{2}{g}+o(1)}$, for any $g$. In \cite{Ci1}, the first author showed that one may take all the $a_k$ to be squares;
and in \cite{Ci2} he showed that there exists an infinite  $B_2[g]$-set $\{ a_k\}$ with  $a_k\ll k^{2+\frac 1g}(\log k)^{\frac 1g+o(1)}$. Here we adapt this latter approach to the set of squares.

\begin{thm}
For any positive integer $g$ there exists an infinite $B_2[g]$
sequence of squares $\{ a_k\}$ such that
$$
a_k\ll k^{2+\frac 1g}(\log k)^{O_g(1)}
$$
\end{thm}

\noindent {\sl Proof}:\ Let $X_1, X_2,\dots$ be an infinite sequence of independent random variables, each of which take values 0 or 1, where
$$
p_b:= \mathbf P(X_b=1) = \ \frac 1{b^{\frac 1{2g+1}} (\log (2+b))^{\beta_g}} .
$$
where $\beta_g>1 $ is a  number we will choose later.
For each selection of random variables we construct a set of integers $\mathcal B=\{ b\geq 1:\ X_b=1\}=\{ b_1<b_2<\dots\}$.  By the central limit theorem we  have $\B(x)\sim c \ x^{1-\frac 1{2g+1}}/(\log x)^{\beta_g}$ with probability $1$ or, equivalently,
$ b_k\sim c' \ (k(\log k)^{\beta_g})^{1+\frac 1{2g}}$.

We will remove from our  sequence of integers $\mathcal B$ any integer $b_0$ such that there exists $n$ for which there are $g+1$ distinct representations of $n$ as the sum of two squares of elements of $\mathcal B$, in which $b_0$ is the very largest element of $\mathcal B$ involved. Let $\mathcal D\subset \mathcal B$ denote the set of such integers $b_0$. Then the set   $\{ c^2: c\in  \mathcal B\setminus \mathcal D \}$ is the desired $B_2[g]$ sequence of squares.

Now, if $b_0\in \mathcal D$ then, by definition,  there exits
$b'_0,b_1,b'_1,\dots b_g,b'_g \in \mathcal B$ with ${b'}_g\le
b_g<\cdots <b_1<b_0$, for which
$$
n=b_0^2+{b'}_0^2=b_1^2+{b'}_1^2=\cdots =b_g^2+{b'}_g^2.
$$
Define  $R(n)=\{ (b,b'),\ b\ge b',\ b^2+b'^2=n\}$, and $r(n)=|R(n)|$. Then the probability that $b_0\in \mathcal D$
because of this particular value of $n$ is
$$
\mathbf E( X_{b_0}X_{{b'}_0} \sum_{\substack{(b_1,{b'}_1),\dots ,(b_g,{b'}_g)\in
R(b_0^2+{b'}_0^2)\\ b_g<\cdots <b_1<b_0}} X_{b_1}X_{{b'}_1}\cdots X_{b_g}X_{{b'}_g} )
.
$$
The $b_j, b_j'$ are all distinct except in the special case that
$n=2b_g^2$ with $b_g=b_g'$. Thus, other than in this special case,
$\mathbf E ( X_{b_0}X_{{b'}_0}  X_{b_1}X_{{b'}_1}\cdots X_{b_g}X_{{b'}_g} )= \prod_{j=0}^g p_{b_j}p_{{b'}_j}\leq (p_{b_0}p_{{b'}_0})^{g+1}$,
since $p_{b_j}p_{{b'}_j}\leq p_{b_0}p_{{b'}_0}$ for all $j$.
This gives a contribution above of
$\leq (p_{b_0}p_{{b'}_0})^{g+1} \binom{r(n)-1}{g}$. The terms with $n=2b_g^2$ similarly contribute
$\leq (p_{b_0}p_{{b'}_0})^{g+1/2} \binom{r(n)-2}{g-1}\leq
p_{b_0}^{2g+1}  r(n)^{g-1}\ll r(n)^{g-1}/b_0'$.
Therefore
$$
\mathbf E(\mathcal D(x)-\mathcal D(x/2)) \ll
\sum_{\substack{{b'}_0\le b_0  \\ x/2<b_0\le x }}
(p_{b_0}p_{{b'}_0})^{g+1} r(b_0^2+{b'}_0^2)^{g}
+  \sum_{\substack{{b'}_0<b<b_0\le x\\ {b'}_0^2+b_0^2=2b^2} } \frac 1 {b_0'}  r(2b^2)^{g-1}  .
$$
For the second sum note that $r(m)\ll m^{o(1)}$ and
that for any $n$ (and in particular for $n=b_0'^2$) we have $\# \{ (y,z),\ n=2z^2-y^2,\ y,z\le x\}\ll (nx)^{o(1)}$, and so  its
total contribution is $\ll x^{o(1)} \sum_{b_0'\le x} 1/b_0=x^{o(1)}$.

For the first term we apply H\"older's inequality with $p=2-\frac 1{g+1}$ and $q=2+\frac 1g$ to   obtain
\begin{eqnarray*}
\leq    \left ( \sum_{\substack{{b'}_0\le b_0  \\ x/2<b_0\le x }} (p_{b_0}p_{{b'}_0})^{2g+1} \right )^{\frac{g+1}{2g+1}}\left (
\sum_{{b'}_0\le b_0\le x}r ^{2g+1}(b_0^2+{b'}_0^2)\right )^{\frac
g{2g+1}} .
\end{eqnarray*}
As $\beta_g>1$, we have
\begin{eqnarray*}
\sum_{\substack{{b'}_0\le b_0  \\ x/2<b_0\le x }} (p_{b_0}p_{{b'}_0})^{2g+1}\ll
\sum_{x/2<b_0\le x} \frac 1{b_0(\log b_0)^{\beta_g(2g+1)} } \sum_{{b'}_0\le b_0} \frac 1{{b'}_0 (\log {b'}_0)^{\beta_g(2g+1)}} \\ \ll \frac 1{ (\log x)^{\beta_g(2g+1)}}  , \\ \text{ and}  \quad  \sum_{{b'}_0\le b_0\le
x}r ^{2g+1}(b_0^2+{b'}_0^2)\le \sum_{n\le 2x^2}r^{2g+2}(n)\ll
x^2(\log x)^{2^{2g+1}-1} ,
\end{eqnarray*}
so that
$$
\mathbf E(\mathcal D(x)-\mathcal D(x/2))  \ll x^{\frac{2g}{2g+1}}(\log x)^{e_g} \ \ \text{where} \ \
e_g:=g\left( \frac{2^{2g+1}-1}{2g+1}\right) -\beta_g(g+1).
$$

Markov inequality's tells us that $\mathbf P\left (\mathcal D(2^j)\ge
j^2\mathbf E(\mathcal D(2^j)-\mathcal D(2^{j-1}))\right )\le  1/j^2$ so that
$\sum_{j\ge 1}\mathbf P\left (\mathcal D(2^j)-\mathcal D(2^{j-1})\gg
j^{2+e_g}(2^j)^{\frac{2g}{2g+1}})\right )<\infty. $ The
Borel-Cantelli lemma then implies that
$$
 \mathcal D(2^j)-\mathcal D(2^{j-1}) \ll
j^{2+e_g}(2^j)^{\frac{2g}{2g+1}} =o(\mathcal B(2^j)-\mathcal B(2^{j-1}) )
$$
with probability $1$, provided $\beta_g>\frac{2^{2g+1}-1}{2g+1} +\frac 2g$. Thus there exists a $B_2[g]$-sequence of the form
$\mathcal A:= \{ a^2:\ a\in \mathcal B\setminus \mathcal D\}$,
where $a_k\ll k^{2+\frac 1g} (\log k)^{\beta_g(1+\frac{1}{2g})}$.

\bigskip

\begin{Cor}
There exists an infinite Sidon sequence of squares $\{a_k\}$ with
$a_k\ll k^3(\log k)^8$.
\end{Cor}

\noindent {\sl Proof}: Take $g=1$ and $\beta=16/3$ in the proof above.

\bigskip

\section{Generalized arithmetic progressions of squares}

A generalized arithmetic progression (GAP) is a set of numbers of the form $\{ x_0+\sum_{i=1}^d j_i x_i: \ 0\leq j_i\leq J_i-1\}$ for some integers $J_1,J_2,\dots ,J_d$ and each $x_i\ne 0$.  We have seen that the questions in this article are closely related to GAPs of squares of integers.  At the start of the article we noted Fermat proved that there are no arithmetic progressions of squares of length 4, and so we may assume each $J_d\leq 3$.
We also saw {Solymosi}'s conjecture \ref{Solymosi} which claims that there are no GAPs of squares with each $J_i=2$ and $d$ sufficiently large. This leaves us just a few cases left to examine:\

We begin by examining arithmetic progressions of length 3 of squares:\  If $x^2,y^2,z^2$ are in arithmetic progression then they satisfy the Diophantine equation $x^2+z^2=2y^2$.  All integer solutions to this equation can be parameterized as
$$
x=r(t^2-2t-1),\ y=r(t^2+1),\ z=r(t^2+2t-1), \ \ \text{where} \ t\in \mathbb Q \ \text{and} \ r \in \mathbb Z.
$$
Therefore the common difference $\Delta$ of this arithmetic progression
is given by $\Delta=z^2-y^2= 4r^2 (t^3-t)$.  Integers which are
a square multiple of numbers of the form $t^3-t,\ t\in \mathbb Q$
are known as {\sl congruent numbers} and have a rich, beautiful history in arithmetic geometry (see Koblitz's delightful book \cite{Kob}). They occur, traditionally, since if a right-angled triangle has rational sides then these can be parameterized as
$s(t^2-1), 2st, s(t^2+1)$ with $s,t\in \mathbb Q$, and so has area $s^2(t^3-t)$ (there is a direct correspondence here since we may take the right-angled triangle to have sides $x+z, z-x, 2y$ which has area $z^2-x^2=2\Delta$). It is a highly non-trivial problem to classify the congruent numbers; indeed this is one of the basic questions of modern arithmetic geometry, see \cite{Kob}.

So can we have a 2-by-3 GAP?  This would require having two different ways to obtain the same congruent number.  The theory of elliptic curves tells us exactly how to do this:\ We begin
with the {\sl elliptic curve}
\begin{equation}\label{ec}
E_\Delta: \ \Delta Y^2 = X^3-X
\end{equation}
and the 3-term arithmetic progressions of rational squares are
in 1-to-1 correspondence with the rational points $(t,1/2r)$ on (\ref{ec}). Now the rational points on an elliptic curve form an abelian group and so if $P=(t,1/2r)$ is a rational point on $E_\Delta$ then there are rational points $2P, 3P,\dots$. This
is all explained in detail in \cite{Kob}. All we need is to note that $2P=(T,1/2R)$ where
$$
T=\frac{(t^2+1)^2}{4(t^3-t)}=\frac{y^2}{\Delta}\ \ \text{and} \ \ R= \frac {8r(t^3-t)^2}{(t^2+1)(t^2+2t-1)(t^2-2t-1)}=\frac {\Delta^2}{2xyz}.
$$
So we have infinitely many 2-by-3 GAPs of squares where the common difference of the 3-term arithmetic progressions is $\Delta$,
for any congruent number $\Delta$.

How about 3-by-3 GAPs of squares? Let us suppose that the common difference in one direction is $\Delta$; having a 3-by-3 GAP is then equivalent to having $y_1^2, y_2^2, y_3^2$  in arithmetic progression. But note that $y_i^2=\Delta T_i= \Delta x(2P_i)$
(where $x(Q)$ denotes the $x$-coordinate of $Q$ on a given elliptic curve). Therefore 3-by-3 GAPs of squares are in 1-to-1 correspondence with the sets of congruent numbers and triples of rational points, $(\Delta; P_1,P_2,P_3): \ P_1,P_2,P_3\in E_\Delta(\mathbb Q)$ for which the $x$-coordinates $x(2P_1), x(2P_2), x(2P_3)$ are in arithmetic progression (other than the triples $-1,0,1$ which do not correspond to squares of interest).

In \cite{BST} it is proved that if there is such an arithmetic progression of rational points then the {\sl rank} of $E_\Delta$ must be at least 2; that is there are at least two points of infinite order in the group of points that are independent.
Bremner became interested in the same issue from a seemingly quite different motivation:

A 3-by-3 {\sl magic square} is a 3-by-3 array of numbers where each row, column and diagonal has the same sum. Solving the linear equations that arise it may be parameterized as
\[ \left( \begin{array}{ccc}
u+v & u-v-\Delta & u+\Delta \\
u-v+\Delta & u & u+v-\Delta \\
u-\Delta & u+v+\Delta & u-v \end{array} \right)\]
The entries of the magic square form the 3-by-3 GAP
$\{ (u-v-\Delta) + j_1 v+ j_2 \Delta:\ 0\leq j_1,j_2\leq 2\}$. Hence
the question of finding a non-trival 3-by-3 magic square with entries from a given set $E$
is {\sl equivalent} to the question of finding a non-trival 3-by-3 GAP with entries
from a given set $E$; in particular when $E$ is the set of  squares.
(This connection is beautifully explained in \cite{Rob}.)

We believe that the existence of non-trivial 3-by-3 GAPs of squares,
and equivalently of non-trivial 3-by-3 magic squares of squares, remain open.

\bigskip

\section{The $abc$-conjecture}

In \cite{BGP} it was shown that the large sieve implies that if there are $\gg \sqrt{k}\log k$ squares amongst $a+b,a+2b,\dots ,a+kb$ then $b\geq e^{\sqrt{k}}$. We wish to obtain an upper bound on $b$ also. We shall do so assuming one of the most important conjectures of arithmetic geometry:

\begin{conj}\label{abc} {\bf (The $abc$-conjecture)} \  If $a+b=c$ where $a,b$ and $c$ are coprime positive integers then
$r(abc)\gg c^{1-o(1)}$ where $r(abc)$ is the product of the distinct primes dividing $abc$.
\end{conj}

Unconditional results on the $abc$-conjecture are from this objective, giving only that $r(abc)\gg (\log c)^{3-o(1)}$, for some $A>0$ (see \cite{SY}). Nonetheless, by considering the strongest feasible version of certain results on linear forms of logarithms, Baker \cite{Bak}  made  a  conjecture which implies the stronger
\begin{equation} \label{bakerabc}
r(abc)\gg c/\exp( (\log c)^{\tau}) ,
\end{equation}
with $\tau=1/2+o(1)$.

\begin{lemma}\label{abcimply} Suppose that $A+t_jB$ is a square for $j=1,2,3,4,5$, where $A,B$ and the $t_j$ are integers
and $(A,B)=1$. Let $T=\max_j |t_j|$. Then (\ref{bakerabc}) implies that $A+B\ll \exp(O(T^{9\tau/(1-\tau)}))$. Moreover if
$B\gg A^{5/6-\epsilon}$ then we may improve this to
$B\ll \exp(O(T^{6\tau/(1-\tau)}))$.
\end{lemma}

\noindent {\sl Proof}:\ There is always a partial fraction decomposition
$$
\frac 1{\prod_{j=1}^5 (x+t_j)} = \sum_{j=1}^5 \frac{e_j}{x+t_j} \ \ \text{where} \ \ e_j = \frac 1{\prod_{i=1,\ i\ne j}^5 (t_i-t_j)} ,
$$
so that $\sum_j e_j t_j^\ell=0$ for $0\leq \ell\leq 3$.
Let $L$ be the smallest integer such that each $E_j:=Le_j$ is an integer.  Define the polynomials
$$
h(x):= \prod_{\substack{1\leq j\leq 5 \\ E_j>0}} (x+t_j)^{E_j} \ \ \text{and} \ \ g(x):= \prod_{\substack{1\leq j\leq 5 \\ E_j<0}} (x+t_j)^{-E_j}, \ \ \text{with} \ \ f(x):=h(x)-g(x).
$$
If $h(x)$ has degree $D$ then the coefficient of $x^{D-i}$ in $f(x)$ is a polynomial in the $\sum_j e_j t_j^\ell$ with $0\leq \ell \leq i$, so we deduce that $f(x)$ has degree $D-4$.
Now let $a=B^Dh(A/B), b=B^Dg(A/B), c=B^4\cdot B^{D-4}f(A/B)$ and
then $a'=a/(a,b), b'=b/(a,b), c'=c/(a,b)$.
Thus $r(a'b'c')\leq r(\prod_{j=1}^5(A+t_jB))|B| |c'/B^4|\leq  \prod_{j=1}^5(A+t_jB)^{1/2}|c'|/B^3$.  Now
$\prod_{j=1}^5(A+t_jB) \ll B^{6-2\epsilon}$ provided
$T=B^{o(1)}$ and $A\ll B^{6/5-\epsilon}$, in which case
$r(a'b'c')\ll |c'|/B^\epsilon$. Then, by (\ref{bakerabc}), we have $(\log c)^{\tau}    \gg \ \log B$.
Now $c=a+b\ll (A+TB)^D$ so that $\log c\ll D\log B$; we deduce that $B\ll \exp(O(D^{\tau/(1-\tau)}))$. Finally note that
$D\ll \max_\ell |E_\ell|  \leq \prod_{1\leq i<j\leq 5,\ i,j\ne \ell} |t_i-t_j| \ll T^6$, and the second result follows.

In case that $A\gg B^{6/5-\epsilon}$ we may replace $t_j$ by $1/t_j$ in our construction of polynomials given above. In that case we get new exponents $e_j^*=e_jt_j^3 \prod_{i=1}^5 t_i$ and therefore $|E_j^*|\leq |t_j|^3E_j$. We now have integers $a^*=\kappa A^dh^*(\frac BA), \ b^*=\kappa  A^dg^*(\frac BA), \ c^*=\kappa  A^4\cdot A^{d-4}f^*(\frac BA)$ where $\kappa:= \prod_j t_j^{|E_j^*|}$ and $d$ is the degree of $h^*$. Thus we have that either $A\ll T^{O(1)}$ or
$A\ll \exp(O(d^{\tau/(1-\tau)}))$ and $d\ll T^9$.
\bigskip

We can apply this directly:\ If there are
$\gg \sqrt{k}$ squares amongst $a+b,a+2b,\dots ,a+kb$
then there must be $i_1<\dots <i_5$ with $i_5<i_1+O(\sqrt{k})$
such that each $a+i_jb$ is a square. Thus by Lemma \ref{abcimply} with $A=a+i_ib,\ B=b,\ t_j=i_j-i_1$, assuming (\ref{bakerabc}) with Baker's $\tau=1/2+o(1)$, we obtain  $a+b\ll \exp(k^{9/2+o(1)})$.  Therefore we may, in future, restrict our attention to the case $k^{1/2} \ll \log(a+b) \ll k^{9/2+o(1)}$.

\bigskip
\noindent{\sl Acknowledgements}:\ Many thanks to  Bjorn Poonen for
his permission to discuss his unpublished work at the end of
section 4.

\end{document}